\documentclass[12pt]{amsart}
\usepackage{amsmath,amssymb,amsthm}

\newcommand{\R}{\mathbb R}
\newcommand{\C}{\mathbb C}

\newtheorem{thm}{Theorem}[section]
\newtheorem{cor}[thm]{Corollary}

\newtheorem{prop}[thm]{Proposition}
\theoremstyle{definition}

\theoremstyle{remark}
\newtheorem{rem}[thm]{Remark}
\newcommand{\ds}{\displaystyle}

\begin{document}

\title[CANONICAL WEIERSTRASS REPRESENTATION OF MINIMAL SURFACES]
{CANONICAL WEIERSTRASS REPRESENTATION OF MINIMAL SURFACES IN EUCLIDEAN SPACE}%

\author{Georgi Ganchev}

\address{Bulgarian Academy of Sciences, Institute of Mathematics and Informatics,
Acad. G. Bonchev Str. bl. 8, 1113 Sofia, Bulgaria}%
\email{ganchev@math.bas.bg}%


\begin{abstract}

Using the fact that any minimal strongly regular surface carries locally canonical
principal parameters, we obtain a canonical representation of these surfaces, which
makes more precise the Weierstrass representation in canonical principal parameters.
This allows us to describe locally the solutions of the natural partial differential
equation of minimal surfaces.

\end{abstract}

\maketitle

\section{Introduction}

In \cite{GM} we proved that any minimal strongly regular surface can be endowed locally
with canonical principal parameters. Using this result, in this note we
prove the following
\vskip 1mm
\noindent
{\bf Theorem 1. (Canonical Weierstrass representation)} {\it Any minimal strongly regular  surface
$\mathcal M: \, \textbf{z}=\textbf{z}(x,y), \; (x,y) \in {\mathcal D} \subset \mathbb{C}$
parameterized with canonical principal parameters has locally a representation of the type
$$\mathcal M: \; \begin{array}{l}
\ds{z_1=\mathfrak{Re}\left( \frac{1}{2} \int_{\bf z_0}^{\bf z}\, \frac{w^2-1}{w'} \, \,dz\right)},\\
[4mm]
\ds{z_2=\mathfrak{Re}\left(-\frac{i}{2} \int_{\bf z_0}^{\bf z} \, \frac{w^2+1}{w'} \, \,dz\right)},\\
[4mm]
\ds{z_3=\mathfrak{Re}\left(-\int_{\bf z_0}^{\bf z} \, \frac{w}{w'} \, \, dz\right)},
\end{array}$$
where
$$w=u(x,y)+iv(x,y), \qquad \mu:=\frac{(u_x^2+u_y^2)}{(u^2+v^2+1)^2}$$
is a holomorphic function satisfying the conditions
$$\mu>0, \quad \mu_x \mu_y \neq 0.$$}

As an application of this theorem we obtain a local description of the solutions of the natural
partial differential equation of minimal surfaces.

\section{Preliminaries}

Let $\mathcal M: \, \textbf{z}=\textbf{z}(x,y), \; (x,y) \in {\mathcal D}$ be a surface
in Euclidean space ${\R}^3$ and $\nabla$ be the flat Levi-Civita connection of the
standard metric in ${\R}^3$. The unit normal vector field to ${\mathcal M}$ is denoted
by $l$ and $E, F, G; \; e, f, g$ stand for the coefficients of the first and the second
fundamental forms, respectively.  All functions are supposed to be in the class
${\mathcal C}^{\infty}$.

The considerations in this note are local.

We suppose that the surface has no umbilical points and the principal lines on
$\mathcal M$ form a parametric net, i.e.
$$F(x,y)=f(x,y)=0, \quad (x,y) \in \mathcal D.$$
Then the principal curvatures $\nu_1, \nu_2$ and the principal geodesic curvatures
(the geodesic curvatures of the principal lines) $\gamma_1, \gamma_2$ are given by
$$\nu_1=\frac{e}{E}, \quad \nu_2=\frac{g}{G}; \qquad
\gamma_1=-\frac{E_y}{2E\sqrt G}, \quad \gamma_2= \frac{G_x}{2G\sqrt E}.$$

We consider the tangential frame field $\{X, Y\}$ determined by
$$X:=\frac{\textbf{z}_x}{\sqrt E}, \qquad Y:=\frac{\textbf{z}_y}{\sqrt G}$$
and assume that the moving frame $XYl$ is always right oriented so that $\nu_1-\nu_2>0$.
The following Frenet type formulas for the frame field $XYl$ are valid
$$\begin{tabular}{ll}
$\begin{array}{llccc}
\nabla_{X} \,X & = &  &\gamma_1 \,Y + \nu_1 \, l,  &\\
[3mm]
\nabla_{X} Y & = -\gamma_1 \, X, & & \\
[3mm]
\nabla_{X} \, l & = - \nu_1 \, X, & & &
\end{array}$ &
\quad
$\begin{array}{llccc}
\nabla_{Y} \,X & = & & \gamma_2 \, Y, &\\
[3mm]
\nabla_{Y} Y & = -\gamma_2 \, X & &  & + \nu_2 \, l,\\
[3mm]
\nabla_{Y}\, l & = &  & - \nu_2 \, Y. &
\end{array}$
\end{tabular}\leqno(2.1)$$
\vskip 2mm
The Codazzi equations are as follows
$$\gamma_1=\frac{Y(\nu_1)}{\nu_1-\nu_2}, \qquad \gamma_2=\frac{X(\nu_2)}{\nu_1-\nu_2},
\leqno(2.2)$$

A surface ${\mathcal M}: \; \textbf{z}=\textbf{z}(x,y), \; (x,y)\in \mathcal D$
parameterized with principal parameters is \emph{strongly regular} \cite{GM} if
$$\gamma_1(x,y)\gamma_2(x,y) \neq 0, \quad (x,y)\in\mathcal D.$$

Because of (2.2)
$$\gamma_1 \gamma_2 \neq 0 \; \iff \; (\nu_1)_y (\nu_2)_x \neq 0.$$

Now let $\mathcal M$ be a minimal strongly regular surface, whose parametric net is
principal. We use the following notations:
$$\nu:=\nu_1>0, \quad \nu_2=-\nu<0, \quad \nu_1-\nu_2=2\nu>0$$
and refer to the function $\nu$ as the normal curvature function.

In \cite {GM} we proved that any minimal strongly regular surface $\mathcal M$ admits
locally canonical principal parameters $(x, y)$ so that the coefficients $E, G$ and
$e, g$ are given by:
$$E=G=\frac{1}{\nu}>0, \quad e=-g=1.\leqno(2.3)$$
\begin{rem} Here we use different from \cite{GM} normalization for $E$ and $G$ in order
to obtain more appropriate form for the canonical representation of minimal surfaces.
\end{rem}
Further we assume, that the minimal strongly regular surface
$\mathcal M: \, \textbf{z}=\textbf{z}(x,y), \; (x,y) \in {\mathcal D}$
is parameterized with canonical principal parameters. Then the principal geodesic
curvatures are given by
$$\gamma_1=(\sqrt{\nu})_y, \quad \gamma_2=-(\sqrt{\nu})_x.\leqno(2.4)$$

The fundamental Bonnet theorem, applied to minimal strongly regular surfaces
parameterized with canonical principal parameters states as follows (cf \cite {GM}):
\vskip 2mm
{\bf Bonnet theorem for minimal surfaces in canonical principal parameters.}
{\it Given a function $\nu (x,y) > 0$
in a neighborhood $\mathcal D$ of $(x_0, y_0)$ with $\nu_x \nu_y \neq 0$,
satisfying the equation
$$\Delta \ln \nu + 2 \nu=0 \leqno(2.5)$$
and an initial right oriented orthonormal frame ${\bf z}_0X_0Y_0l_0$.

Then there exists a unique minimal strongly regular surface
${\mathcal M}: \; \textbf{z}=\textbf{z}(x, y), \; (x, y) \in \mathcal D_0 \;
((x_0, y_0) \in \mathcal D_0 \subset \mathcal D)$, such that

$(i)$ \; \; $(x, y)$ are canonical principal parameters;

$(ii)$ \, \, ${\bf z}(x_0, y_0)={\bf z}_0, \; X(x_0, y_0)=X_0, \;
Y(x_0, y_0)=Y_0, \; l(x_0, y_0)=l_0$;

$(iii)$ \, the invariants of ${\mathcal M}$ are the following functions
$$\nu_1=\nu, \quad \nu_2= -\nu, \quad \gamma_1=(\sqrt {\nu})_y,  \quad
\gamma_2= -(\sqrt{\nu})_x.$$}

Further we refer to (2.5) as the natural partial differential equation of
minimal surfaces.

The above statement gives a one-to-one correspondence between minimal strongly regular
surfaces (considered up to a motion) and the solutions of the natural
partial differential equation, satisfying the conditions
$$\nu > 0, \quad \nu_x\nu_y\neq 0. \leqno(2.6)$$

According to \cite{GM} the sign of the function $\nu_x \nu_y$ divides the class of minimal
strongly regular surfaces into two geometric subclasses (invariant with respect to motions
and changes of parameters). Any reflection of ${\R}^3$ with respect to a plane transforms
each of these subclasses onto the other.

\section{Canonical representation of minimal strongly regular surfaces}

Let $\mathcal M: \, \textbf{z}=\textbf{z}(x,y), \; (x,y) \in {\mathcal D}$ be a minimal
strongly regular surface parameterized with canonical principal parameters. In view of
(2.3) and (2.4) formulas (2.1) become
$$\begin{array}{llll}
\nabla_{X} \, X & = &\;\;\;(\sqrt{\nu})_y  \, Y & + \, \nu \, l, \\
[2mm]
\nabla_{X} \, Y & = -(\sqrt{\nu})_y \, X, \\
[2mm]
\nabla_{X} \, l &= - \nu \, X;\\
[4mm]
\nabla_{Y} \, X & = & -(\sqrt{\nu})_x \, Y, \\
[2mm]
\nabla_{Y} \, Y & = \;\;\;(\sqrt{\nu})_x \, X & & - \, \nu \, l,\\
[2mm]
\nabla_{Y} \, l & = & \;\;\; \nu \, Y
\end{array}
\leqno(3.1)$$
and the integrability conditions for (3.1) reduce to (2.5).

Equalities (3.1) imply the following formulas for the Gauss map $l=l(x,y)$:
$$\begin{array}{ll}
l_{xx} &= \;\;\; \ds{\frac{\nu_x}{2\nu} \, l_x \, - \,
\frac{\nu_y}{2 \nu}\, l_y \, - \, \nu \, l,}\\
[4mm]
l_{xy} &= \;\;\; \ds{\frac{\nu_y}{2\nu}\, l_x \, + \, \frac{\nu_x}{2 \nu}\,l_y,}\\
[4mm]
l_{yy} &= - \ds{\frac{\nu_x}{2\nu}\,l_x \, + \, \frac{\nu_y}{2 \nu}\, l_y
\, - \,\nu \, l}
\end{array} \leqno(3.2)$$
and the unit normal vector function $l(x,y)$ satisfies the differential equation:
$$\Delta l+2l=0.$$

The next statement makes precise the relation between the properties of the Gauss map of
a minimal surface and the canonical principal parameters.
\begin{prop}\label{P:3.1}
Let $\mathcal M: \textbf{z}=\textbf{z}(x,y), \; (x,y) \in {\mathcal D}$ be a minimal
strongly regular surface parameterized with canonical principal parameters. Then the
Gauss map $l=l(x,y), \; (x,y)\in {\mathcal D}; \; l^2=1$ has the following properties:
$$l_x^2=l_y^2=\nu >0, \quad l_x \,l_y=0,\quad \nu_x\nu_y\neq 0.\leqno(3.3)$$

Conversely, if a unit vector function $l(x,y)$ has the properties $(3.3)$, then there
exists locally a unique (up to a motion) minimal strongly regular surface
$\mathcal M :\; \textbf{z}=\textbf{z}\,(x,y)$ determined by
$$\textbf{z}_x=-\frac{1}{\nu}\,l_x, \quad \textbf{z}_y=\frac{1}{\nu}\,l_y,
\leqno(3.4)$$
so that $(x,y)$ are canonical principal parameters and $\nu(x,y)$ is the normal curvature
function of $\mathcal M$.
\end{prop}
\emph{Proof:} The equalities  $l_x=-\nu \, \textbf{z}_x, \; l_y=\nu \, \textbf{z}_y$,
and (2.3) imply (3.3).

For the inverse, it follows immediately that (3.3) implies (3.2). Therefore the system
(3.4) is integrable and determines locally a surface
$\mathcal M:\;\textbf{z}=\textbf{z}(x,y)$.

Since (3.4) and (3.2) imply (3.1), it follows that $\mathcal M$ is a minimal strongly
regular surface parameterized with canonical principal parameters, whose normal
curvature function is $\nu=\,l_x^2=\,l_y^2 >0$.

Furthermore, it follows that the function $\nu=\,l_x^2=\,l_y^2$ satisfies the equation
(2.5). $\qed$
\vskip 2mm
Thus, any minimal strongly regular surface locally is determined by the system (3.4), where
the unit vector function $l=l(x,y)$ satisfies the conditions (3.3). The so obtained minimal
surface is parameterized with canonical principal parameters.
\vskip 2mm
Let $S^2(1):\xi^2 + \,\eta^2 + \,\zeta^2=1$ be the unit sphere centered at the origin
$O$ and \, $l(\xi, \eta, \zeta), \\ \zeta \neq 1$ be the position vector of an arbitrary
point on $S^2(1)$, different from the pole $(0,0,1)$. The standard conformal parametrization
of $S^2$ generated by the stereographic map
$$l: \quad \begin{array}{l}
\displaystyle{\xi=\frac{2x}{x^2+y^2+1}},\\
[4mm]
\displaystyle{\eta=\frac{2y}{x^2+y^2+1}},\\
[4mm]
\displaystyle{\zeta=\frac{x^2+y^2-1}{x^2+y^2+1}},
\end{array}\qquad (x,y) \in {\R}^2 $$
satisfies the equalities
$$l_x^2=l_y^2=\frac{4}{(x^2+y^2+1)^2}, \quad l_x\,l_y=0$$
and the unit vector function $l=l(x,y)$ generates the Enneper minimal surface parameterized
with canonical principal parameters.

Now we make more precise the Weierstrass representation of minimal strongly regular surfaces
in canonical principal parameters.
\vskip 4mm
\noindent
{\bf Proof of Theorem 1}
\vskip 4mm
Let $\mathcal M: \textbf{z}=\textbf{z}(x,y), (x,y) \in \mathcal D$ be
a minimal strongly regular surface parameterized with canonical principal parameters.
Since the Gauss map $l=l(x,y)$ of $\mathcal M$ is conformal and the orthogonal frame field
$l_u \, l_v \, l$ is left oriented, then the vector function $l$ is given locally by
the equalities
$$l: \quad \begin{array}{l}
\displaystyle{\xi=\frac{2u(x,y)}{u^2(x,y)+v^2(x,y)+1}},\\
[4mm]
\displaystyle{\eta=\frac{2v(x,y)}{u^2(x,y)+v^2(x,y)+1}},\\
[4mm]
\displaystyle{\zeta=\frac{u^2(x,u)+v^2(x,y)-1}{u^2(x,y)+v^2(x,y)+1}},
\end{array}\leqno(3.5)$$
where
$$w: \quad
\begin{array}{l}
u=u(x,y),\\
[2mm]
v=v(x,y),\end{array}
\qquad \qquad
\begin{array}{l}
u_x=\;\;\;v_y,\\
[2mm]
u_y=-v_x
\end{array}
$$
is a holomorphic function in ${\C}$.

We denote
$$\mu :=\frac{(u_x^2+u_y^2)}{(u^2+v^2+1)^2}.$$

If $\nu=\nu(x,y)$ is the normal curvature function of $\mathcal M$, then
the vector function $\textbf{z}=\textbf{z}(x,y)$ satisfies the system
$$\begin{array}{l}
\ds{\textbf{z}_x=-\frac{1}{\nu}\,l_x=-\frac{1}{\nu}\,(u_x\,l_u-u_y \,l_v)},\\
[4mm]
\ds{\textbf{z}_y=\;\;\;\frac{1}{\nu}\,l_y=\;\;\;\frac{1}{\nu}\,(u_y\,l_u+u_x\,l_v).}
\end{array}\leqno(3.6)$$

It follows from (3.6) that $\ds{\nu=\frac{4(u_x^2+u_y^2)}{(u^2+v^2+1)^2}}=4\mu$.
Hence the holomorphic function $w$ satisfies the conditions
$$\mu > 0, \quad \mu_x \mu_y \neq 0.$$

Denoting by $w'=\ds{\frac{dw}{dz}=\frac{\partial w}{\partial z}=u_x-iu_y}$, we find from (3.6)
$${\bf z}_x-i{\bf z}_y=-\frac{1}{w'} \, \frac{(u^2+v^2+1)^2}{4} \, (l_u+il_v).\leqno(3.7)$$
Then (3.7) in view of (3.5) can be written in the form:
$$\begin{array}{l}
\ds{(z_1)_x-i(z_1)_y=\;\;\,\frac{1}{2} \, \frac{w^2-1}{w'}},\\
[4mm]
\ds{(z_2)_x-i(z_2)_y=-\frac{i}{2} \, \frac{w^2+1}{w'}},\\
[4mm]
\ds{(z_3)_x-i(z_3)_y=-\frac{w}{w'}},
\end{array}$$
which proves the assertion. $\hfill{\qed}$
\vskip 2mm
As an application we obtain a corollary for the solutions of the natural
partial differential equation of minimal surfaces.

\begin{cor} Any solution $\nu$ of the natural partial differential
equation $(2.5)$ satisfying the condition $(2.6)$ locally is given by the formula
$$\nu(x,y)=\frac{4(u_x^2+u_y^2)}{(u^2+v^2+1)^2} \leqno(3.8)$$
where $w=u(x,y)+iv(x,y)$ is a holomorphic function in $\mathbb C$.
\end{cor}
{\it Proof:} Let $\nu(x,y)$ be a solution to (2.5) satisfying the condition
(2.6). Then the function $\nu$ generates locally a minimal strongly regular surface
$\mathcal M$ (unique up to a motion). According to Theorem 1 it follows that the normal
curvature function $\nu$ of $\mathcal M$ has locally the form (3.8). \qed

It is a direct verification that any function $\nu$ given by (3.8), where
$w=u+iv \; (u_x^2+u_y^2>0)$ is a holomorphic function, satisfies (2.5).
\vskip 2mm
\begin{rem} The canonical Weierstrass representation is based on the Gauss map of the
Enneper surface ($w=z$). It is clear that choosing the Gauss map of any other minimal
strongly regular surface $\mathcal M$, we shall obtain its corresponding representation.
This remark is also valid for the form (3.8) of the solutions of the natural partial
differential equation (2.5).
\end{rem}

\end{document}